\documentclass[12pt]{article}
\usepackage{amssymb}
\usepackage{amsmath,graphicx,amsfonts,epsfig,eucal}
\usepackage{amsfonts}

\newcommand{\ds}{\displaystyle}
\newcommand{\pa}{\partial}
\newcommand{\rr}{{\rm{I \! R}}}
\newtheorem{theorem}{Theorem}

\newtheorem{proposition}[theorem]{Proposition}

\input{tcilatex}

\begin{document}

\title{The analysis of the stochastic stability for an economic game}
\author{A. L. Ciurdariu$^{a}$, M. Neam\c{t}u$^{b}$, A. Sandru$^{c}$, D. Opri\c s$^{d}$}
\date{}
\maketitle

\begin{tabular}{cccccccc}
\scriptsize{$^{a}$ Department of Mathematics, Politehnica University of Timi\c{s}oara}\\
\scriptsize{P-\c{t}a. Victoriei, nr 2,
300004, Timi\c{s}oara, Romania, e-mail: cloredana43@yahoo.com,}\\
%\scriptsize{E-mail: cloredana43@yahoo.com}\\
\scriptsize{$^{b}$Department of Economic Informatics, Mathematics and Statistics,}\\
\scriptsize{Faculty of Economics, West University of Timi\c soara,}\\
\scriptsize{str. Pestalozzi, nr. 16A, 300115, Timi\c soara, Romania, e-mail:mihaela.neamtu@fse.uvt.ro}\\
%\scriptsize{E-mail:mihaela.neamtu@fse.uvt.ro,}\\
\scriptsize{$^{c}$ Faculty of Exact Sciences, Aurel Vlaicu
University of Arad},\\
\scriptsize{str. Elena Dragoi nr. 2, 310330, Arad, e-mail: sandruandrea@gmail.com,}\\
%\scriptsize{E-mail: sandruandrea@gmail.com,}\\
\scriptsize{$^{d}$ Department of Applied
Mathematics, Faculty of Mathematics,}\\
\scriptsize{West University of Timi\c soara, Bd. V. Parvan, nr. 4, 300223, Timi\c soara, Romania,}\\
 \scriptsize{E-mail: opris@math.uvt.ro}\\

\end{tabular}

\begin{abstract} In this paper we investigate a
stochastic model for an economic game. To describe this model we
have used a Wiener process, as the noise has a stabilization effect.
The dynamics are studied in terms of stochastic stability in the
stationary state, by constructing the Lyapunov exponent, depending
on the parameters that describe the model. Also, the Lyapunov
function is determined in order to analyze the mean square
stability. The numerical simulation that we did justifies the
theoretical results.
\end{abstract}

\medskip

%\newline
{\small \textit {Mathematics Subject Classification}:
34D08,60H10,91B70}

%\newline

{\small \textit {Keywords}: stochastic dynamics in economic games,
economic games, stochastic stability, Lyapunov exponent, Euler
scheme}

\section{Introduction.}

\qquad Stochastic modeling plays an important role in many branches
of science. In many practical situations, perturbations are
expressed in terms of white noise, modeled by brownian motion. The
behavior of a deterministic dynamical system which is disturbed by
noise may be modeled by a stochastic differential equation (SDE),
\cite{KP}. The stochastic stability has been introduced by Bertram
and Sarachik and is characterized by the negativeness of Lyapunov
exponents. In general, it is not possible to determine this
exponents explicitly. Many numerical approaches have been proposed,
which generally used the simulation of the stochastic trajectories.
In the present paper, we study a stochastic dynamical system that is
used in economy, in describing a Counot duopoly game.

In 1838, Cournot introduced the first formal theory of oligopoly,
which treated the case of naive expectations, where each player
assumes the last values taken by the competitors without estimation
of their future reactions \cite{Cournot}. Recently, a lot of
articles have shown that the Cournot model may lead to a cyclic or
chaotic behavior \cite{Bundau}, \cite{Bischi}, \cite{Mircea},
\cite{Puu1}, \cite{Puu2}, \cite{Puu3}. Also, in \cite{Rosser},
Rosser reviews the development of the theory of complex oligopoly
dynamics.

In the present paper we have studied a stochastic Cournot economic
game. In Section 2 we present the Lyapunov exponent and stability in
stochastic 2d dynamical structures. Section 3 describes the Lyapunov
function method for the stochastic stability analysis. Section 4
studies the Lyapunov exponent for an economic game with stochastic
dynamics. The Lyapunov function method for the stochastic game is
given in Section 5. Some numerical simulations are done in Section
6. Finally, Section 7 draws some conclusions.

\section{The Lyapunov exponent and stability in stochastic 2d dynamical structures.}

\qquad Let $(\Omega , {\cal F}, {\cal P})$ be a probability space
\cite{KP}. It is
 assumed that the $\sigma -$algebra ${\cal F}$ is a filtration that
 is, ${\cal F}$ is generated by a family of $\sigma -$algebra ${\cal F}_t(t\geq
 0)$ such that

$${\cal F}_s\subset {\cal F}_t\subset {\cal F}, \quad \forall s\leq t, s,t\in
I,$$ where $I=[0, T]$, $T\in (0, \infty)$.

Let $\{x(t)=(x_1(t), x_2(t))\}_{t\geq 0}$ be a stochastic process.
The system of It$\hat{o}$ equations:
\begin{equation}\label{1}
dx_i(t,\omega)=f_i(t,x(t,\omega))dt+g_i(x(t,\omega))dw(t,\omega),
i=1,2,
\end{equation} with the initial condition $x(0)=x_0$ is written as:
\begin{equation}\label{2}
x_i(t,\omega)=x_{i0}(\omega)+\int_0^tf_i(x(s,\omega))ds+\int_0^tg_i(x(s,\omega))dw(s,\omega),
i=1,2,
\end{equation} for almost all $\omega\in\Omega $ and for each $t>0$,
where $f_i(x)$ is a drift function, $g_i(x)$ is a diffusion
function, $\int_0^tf_i(x(s))ds$, $i=1,2$ is a Riemann integral and
$\int_0^tg_i(x(s))dw(s)$ is an It$\hat{o}$ integral. It is assumed
that $f_i$ and $g_i$, $i=1,2$ satisfy the conditions of existence of
solution for this SDE with initial condition $x(0)=a_0\in\rr^n$.

Let $x_0=(x_{10}, x_{20})\in\rr^2$ be a solution of the system:
\begin{equation}\label{3}
f_i(x_0)=0, i=1,2.
\end{equation}

The functions $g_i, i=1,2$ are chosen so that:
\begin{equation*}
g_i(x_0)=0, i=1,2.
\end{equation*}

In what follows, we consider:
\begin{equation*}\label{4}
g_i(x)=\sum_{j=1}^{2}b_{ij}(x_j-x_{0j}), i=1,2,
\end{equation*} where $b_{ij}\in\rr, i,j=1,2.$

The linearized system of (\ref{2}) in $x_0$, is given by:
\begin{equation*}\label{5}
X(t)=\int_0^tAX(s)ds+\int_0^tBX(s)dw(s),
\end{equation*} where
\begin{equation*}\label{6}
X(t)=\left (\begin{array}{c} u_1(t,\omega )\\ u_2(t,\omega
)\end{array}\right ), A=\left (\begin{array}{cc} a_{11} & a_{12}\\
a_{21} & a_{22}\end{array}\right ), B=\left (\begin{array}{cc} b_{11} & b_{12}\\
b_{21} & b_{22}\end{array}\right ),
\end{equation*}
\begin{equation*}\label{7}
a_{ij}=\ds\frac{\pa f_i}{\pa x_j}|_{x_0}, b_{ij}=\ds\frac{\pa
g_i}{\pa x_j}|_{x_0}.
\end{equation*}

The Oseledec multiplicative ergodic theorem \cite{Ose} asserts the
existence of 2 non-random Lyapunov exponents
$\lambda_2\leq\lambda_1=\lambda$. The top Lyapunov exponent is given
by:
\begin{equation*}\label{8}
\lambda =\lim_{t\to\infty} \sup\log\sqrt{u_1(t)^2+u_2(t)^2}.
\end{equation*}

Applying the change to polar coordinates:
\begin{equation*}
x(t)=r(t)cos \theta (t), y(t)=r(t) sin\theta (t)
\end{equation*} by writing the It$\hat{o}$ formula for
\begin{equation*}
h_1(u_1,u_2)=\ds\frac{1}{2}\log (u_1^2+u_2^2)=\log (r),
h_2(u_1,u_2)=arctg (\ds\frac{u_2}{u_1})=\theta .
\end{equation*} we get:

\begin{proposition}\label{P1} \cite{KP}.  The formulas
\begin{equation}\label{11}
\log \left (\ds\frac{r(t)}{r(0)}\right
)\!=\!\int_0^tq_1(\theta(s))\!+\!\ds\frac{1}{2}(q_4(\theta(s))^2\!-\!q_2(\theta(s))^2)ds\!+\!\int_0^tq_2(\theta(s))dw(s),
\end{equation}
\begin{equation}\label{12}
\theta(t)\!=\!\theta(0)+\int_0^tq_3(\theta(s))\!-\!q_2(\theta(s)q_4(\theta(s))ds\!+\!\int_0^tq_4(\theta(s))dw(s),
\end{equation} hold,
 where
\begin{equation}\label{13}
\begin{array}{llll}
q_1(\theta)=a_{11}cos^2(\theta)+(a_{12}+a_{21})cos\theta\sin\theta+a_{22}sin^2\theta,\\
q_2(\theta)=b_{11}cos^2(\theta)+(b_{12}+b_{21})cos\theta\sin\theta+b_{22}sin^2\theta,\\
q_3(\theta)=a_{21}cos^2(\theta)+(a_{22}-a_{11})cos\theta\sin\theta-a_{12}sin^2\theta,\\
q_4(\theta)=b_{21}cos^2(\theta)+(b_{22}-b_{11})cos\theta\sin\theta-b_{12}sin^2\theta.\\
\end{array}
\end{equation}
\end{proposition}
As the expectation of the It$\hat{o}$ stochastic integral is null
$$E\int_0^tq_2(\theta (s))dw(s)=0,$$ the Lyapunov exponent is given
by:
$$\lambda\!=\!\lim_{t\to\infty}\ds\frac{1}{t}\log\left (\ds\frac{r(t)}{r(0)}\right )\!=\!
\lim_{t\to\infty}\ds\frac{1}{t}E\int_0^t(q_1(\theta(s))\!+\!\ds\frac{1}{2}(q_4(\theta(s)))^2\!-\!q_2(\theta(s)))ds.$$
Applying the Oseledec theorem, if $r(t)$ is ergodic, we get:
\begin{equation*}\label{14}
\lambda
=\int_0^t(q_1(\theta)+\ds\frac{1}{2}(q_3(\theta)^2-q_2(\theta)))p(\theta)d\theta
,
\end{equation*}where $p(\theta)$ is the density of probability of
the process $\theta$.

An approximation of this density is calculated by solving the
Fokker-Planck equation.

The Fokker-Planck (FPE) equation associated with equation (\ref{12})
for $p=p(t,\theta)$ is
\begin{equation}\label{15}
\ds\frac{\pa p}{\pa t}+\ds\frac{\pa }{\pa\theta
}((q_3(\theta)-q_2(\theta)q_4(\theta))p)-\ds\frac{1}{2}\ds\frac{\pa
^2}{\pa \theta^2}(q_4(\theta)^2p)=0.
\end{equation}

From (\ref{15}), it results that the solution $p(\theta)$ of the FPE
is a solution of the following first order equation:
\begin{equation}\label{16}
(-q_3(\theta)\!+\!q_1(\theta)q_4(\theta)\!+
\!q_2(\theta)g_5(\theta))p(\theta)\!+\!\ds\frac{1}{2}q_4(\theta)^2p'(\theta)\!=\!p_0,
\end{equation} where $p'(\theta )=\ds\frac{dp}{d\theta}$ and
\begin{equation*}\label{17}
q_5(\theta)=-(b_{12}+b_{21})sin 2\theta -(b_{22}-b_{11})cos 2\theta
.
\end{equation*}

\begin{proposition}\label{prop2}\cite{KP}. If $q_4(\theta)\neq 0$, the
solution of the equation (\ref{16}) is given by:
\begin{equation*}\label{18}
p(\theta)=\ds\frac{k}{D(\theta)q_4(\theta)^2}\left (1+\eta
\int_0^\theta D(u)du\right )
\end{equation*} where $k$ is determined by the normality condition
\begin{equation*}\label{18}
\int_0^{2\pi }p(\theta )d\theta =1
\end{equation*} and
\begin{equation*}\label{19}
\eta =\ds\frac{D(2\pi )-1}{\int_0^{2\pi }D(u)du}.
\end{equation*}
The function $D$ is given by:
\begin{equation*}\label{21}
D(\theta )=\exp
(-2\int_0^\theta\ds\frac{q_3(u)-q_2(u)q_4(u)-q_4(u)q_5(u)}{q_4(u)^2}du)
\end{equation*}
\end{proposition}

A numerical solution of the phase distribution could be given by a
simple backward difference scheme.

We consider $N\in\rr_+$, $h=\ds\frac{\pi}{N}$ and
\begin{equation*}\label{13}
\begin{array}{llll}
q_1(i)=a_{11}\cos^2(ih)+(a_{12}+a_{21})\cos(ih)\sin(ih)+a_{22}\sin^2(ih),\\
q_2(i)=b_{11}\cos^2(ih)+(b_{12}+b_{21})\cos(ih)\sin(ih)+b_{22}\sin^2(ih),\\
q_3(i)=a_{21}\cos^2(ih)+(a_{22}-a_{11})\cos(ih)\sin(ih)-a_{12}sin^2(ih),\\
q_4(i)=b_{21}\cos^2(ih)+(b_{22}-b_{11})\cos(ih)\sin(ih)-b_{12}\sin^2(ih),\\
q_5(i)=-(b_{12}+b_{21})\sin(2ih)-(b_{22}-b_{11})cos(2ih), i=0,...,N\\
\end{array}
\end{equation*}

The function $p(i), i=0,...,N$ is given by the following relations:
$$p(i)=(p(0)+\ds\frac{q_4(i)^2p(i-1)}{2h})F(i)$$ where
$$F(i)=\ds\frac{2h}{2h(-q_3(i)+q_2(i)q_4(i)+q_4(i)q_5(i))+q_4(i)^2}.$$
The Lyapunov exponent is $\lambda =\lambda (N)$, where
$$\lambda (N)=\sum_{i=0}^{N}(q_1(i)+\ds\frac{1}{2}(q_4(i)^2-q_2(i)^2))p(i)h.$$

From Proposition \ref{prop2} we obtain:

\begin{proposition}\label{prop3} If the matrix B is given by:
$$b_{11}=\alpha, b_{12}=-\beta, b_{21}=\beta, b_{22}=\alpha $$
then
\begin{equation*}
\begin{split}
p(\theta)&=\ds\frac{k}{\beta^2}\exp\{\ds\frac{1}{\beta^2}((a_{21}-a_{12}-\alpha\beta)\theta+\ds\frac{1}{2}(a_{11}-a_{22})\cos
2\theta+ \ds\frac{1}{2}(a_{21}-\\
&-a_{12})\sin2\theta)\}
\end{split}
\end{equation*}
\begin{equation*}
k\!=\!\ds\frac{\beta^2}{\int_0^{2\pi}\exp\{\ds\frac{1}{\beta^2}((a_{21}\!\!-\!\!a_{12}\!\!-\!\!\alpha\beta)\theta\!+\!\ds\frac{1}{2}(a_{11}\!-\!a_{22})\cos
2\theta\!\!+\!\!
\ds\frac{1}{2}(a_{21}\!\!-\!\!a_{12})\sin2\theta)d\theta}
\end{equation*}
\begin{equation*}
\lambda=\ds\frac{1}{2}(a_{11}+a_{22}+\beta^2-\alpha^2)+\ds\frac{1}{2}(a_{11}-a_{22})c_2+\ds\frac{1}{2}(a_{21}+a_{12})s_2,
\end{equation*} where
$$c_2=\int_0^{2\pi}cos(2\theta)p(\theta)d\theta, \quad s_2=\int_0^{2\pi} sin(2\theta)p(\theta)d\theta .$$
\end{proposition}

\section{The Lyapunov function method for the stochastic stability analysis.}
\qquad We consider the system of stochastic equations, SDE, given
by:
\begin{equation}\label{31f}
dx_i(t)=f_i(x(t))dt+g_i(x(t))dB^i(t), \quad i=1,2,
\end{equation} where $x(t)=x(t,\omega )$ and $B^1(t)$, $B^2(t)$ are Wiener processes.

Let $V:D=(0, \infty)\times \rr^2\to\rr$ be a continuous function
with respect to the first variable and a $C^2$ class function with
respect to the other variables.

Let:
\begin{equation}\label{32f}
LV(t,x)=\ds\frac{\pa V(t,x)}{\pa t}+\sum_{i=1}^2f_i(x)\ds\frac{\pa
V(t,x)}{\pa x_i}+\ds\frac{1}{2}\sum_{i=1}^2g_i(x)g_j(x)\ds\frac{\pa
^2V(t,x)}{\pa x_i\pa x_j}
\end{equation} be a differential operator.

We suppose that $x_e=0$ is the stationary state of (\ref{31f}), that
means:
\begin{equation*}
f_i(0)=g_{i\alpha }(0)=0, \quad i=1,2, \alpha =1,2.
\end{equation*}

The theorem which gives us conditions for the stability of the
trivial solution $x_e=0$ in the terms of the Lyapunov function is:

\begin{theorem}\label{th1} \cite{Schu}
Under the above conditions, if there is a function $V:D\to\rr $ and
two continuous functions $u,v:\rr _+\to\rr _+$ and $k>0$ so that for
$||x||<k$ the relation:
\begin{equation*}
u(||x||)<V(t,x)<v(||x||)
\end{equation*} holds, then:

(i) if $LV(t,x)\leq 0$, $x\in (0,k)$, then solution of (\ref{31f})
$x_e=0$ is stable in probability;

(ii) if there is a continuous function $c:\rr _+\to \rr _+$ so that
$LV(t,x)\leq -c(||x||)$ then solution $x_e=0$ of (\ref{31f}) is
asymptotically stable.
\end{theorem}

Let $V:D=(0, \infty)\times \rr^2\to\rr$ be a continuous function
with respect to the first variable and a $C^2$ class function with
respect to the other variables.

The theorem that gives us the exponential p-stability condition of
the trivial solution (\ref{35f}) is:
\begin{theorem}\label{th2} \cite{Schu}
If function $V$ satisfies the following inequalities:
\begin{equation*}
\begin{split}
& k_1||u||^p\leq V(t,x)\leq k_2||u||^p\\
& LV(t,u)\leq -k_3||u||^p, \quad k_i>0, p>0, i=1,2,
\end{split}
\end{equation*} then the trivial solution of (\ref{35f}) is
exponential p-stable for $t\geq 0$.
\end{theorem}

For the concrete problems the following theorem is used:

\begin{theorem}\label{th3} \cite{Hasm}
If function $V$ satisfies the following inequality:

(i) $LV(u)\leq 0$, then the trivial solution is stable in
probability;

(ii) $LV(u)\leq -c(||u||)$, where $c:\rr _+\to\rr_+$ is a continuous
function, then the trivial solution is asymptotically stable;

(iii) $LV(u)\leq -q^TQq$, where $Q$ is a symmetric matrix positively
defined, then the trivial solution is mean square stable.
\end{theorem}

In general, the functions $f_i$, $g_{i\alpha }$, $i=1,2$, $\alpha
=1,2$ are nonlinear functions and the above theorem is difficult to
use. Therefore, the linearization method of system (\ref{31f}), in
the neighborhood of the equilibrium point is used.

The linearized stochastic differential system SDEL of (\ref{31f}) is
given by:

\begin{equation}\label{35f}
\begin{split}
& du_1(t)=(a_{11}u_1(t)+a_{12}u_2(t))dt+(b_{11}u_1(t)+b_{12}u_2(t))dB^1(t)\\
&
du_2(t)=(a_{21}u_1(t)+a_{22}u_2(t))dt+(b_{21}u_1(t)+b_{22}u_2(t))dB^2(t).
\end{split}
\end{equation}

For (\ref{35f}) expression LV is given by:
\begin{equation}\label{36f}
\begin{split}
& LV=(a_{11}u_1+a_{12}u_2)\ds\frac{\pa V}{\pa
u_1}+(a_{21}u_1+a_{22}u_2)\ds\frac{\pa V}{\pa u_2}+\\
& \ds\frac{1}{2}[(b_{11}u_1+b_{12}u_2)^2\ds\frac{\pa ^2 V}{\pa
u_1^2}+(b_{21}u_1+b_{22}u_2)^2\ds\frac{\pa ^2V}{\pa u_2^2}]
\end{split}
\end{equation}

\section{The Lyapunov exponent for an economic game with stochastic dynamics.}
\qquad Two firms enter the market with a homogenous consumption
product. The elements which describe the model are: the quantities
which enter the market
from the two firms $x_{i}\geq 0,$ $i=\overline{1,2};$  the inverse demand function $p:\mathbb{R}%
_{+}\rightarrow \mathbb{R}_{+}$ ($p$ is a derivable function with $%
p^{\prime }\left( x\right) <0,\underset{x\rightarrow
a_{1}}{lim}p\left( x\right) =0,$ $\underset{x\rightarrow 0}{lim}
p\left( x\right) =b_{1},\left( a_{1}\in
\overline{\mathbb{R}},b_{1}\in \overline{\mathbb{R}}\right) $;  the cost functions $C_{i}:%
\mathbb{R}_{+}\rightarrow \mathbb{R}_{+}$ ( $C_{i}$ are derivable
functions with $C_{i}^{\prime }\left( x_{i}\right) >0,$
$C_{i}^{\prime \prime }\geq 0,$ $i=\overline{1,2}$ ).

In our study we consider $p(x)=\ds\frac{1}{x}, x>0$ and
$C_i(x_i)=c_ix_i+d_i, i=1,2.$

The mathematical model of the stochastic dynamic economic game is
described by the stochastic system of equations:
\begin{equation}\label{31}
\begin{split}
x_1(t)\!=\!x_1\!(0)\!\!+\!k_1\!\int_0^t\!(\ds\frac{x_2(s)}{(x_1(s)\!+\!x_2(s))^2}\!-\!c_1)ds\!+\!\int_0^t(b_{11}x_1(s)\!+\!b_{12}x_2(s)\!+\!\gamma_1)dw(s)\\
x_2(t)\!=\!x_2\!(0)\!\!+\!k_2\!\int_0^t\!(\ds\frac{x_1(s)}{(x_1(s)\!+\!x_2(s))^2}\!-\!c_2)ds\!+\!\int_0^t(b_{21}x_1(s)\!+\!b_{22}x_2(s)\!+\!\gamma_2)dw(s)
\end{split}
\end{equation} where $b_{ij}\in\rr$, $i,j=1,2$, $k_1>0, k_2>0$, $x_i(t)=x_i(t,\omega
)$, $i=1,2$.
\begin{equation*}\label{32}
\gamma_1=-\ds\frac{b_{11}c_2+b_{12}c_1}{(c_1+c_2)^2},
\gamma_2=-\ds\frac{b_{21}c_2+b_{22}c_1}{(c_1+c_2)^2}.
\end{equation*}

For $b_{ij}=0$, $i,j=1,2$ model (\ref{31}) is reduced to the
classical model of the economic game \cite{Bundau}, \cite{Mircea}.

The system of stochastic equations (\ref{31}), has the form
(\ref{2}) from section 2, where:
\begin{equation*}\label{33}
\begin{split}
f_1(x_1,x_2)=\ds\frac{x_2}{(x_1+x_2)^2}-c_1,
g_1(x_1,x_2)=b_{11}x_1+b_{12}x_2+\gamma_1,\\
f_2(x_1,x_2)=\ds\frac{x_1}{(x_1+x_2)^2}-c_2,
g_2(x_1,x_2)=b_{21}x_1+b_{22}x_2+\gamma_2.
\end{split}
\end{equation*}

Applying the results from section 2, we have:

\begin{proposition}\label{prop4}
(i) The stationary state of (SDE)
(\ref{31}) is given by:
\begin{equation*}\label{34}
x_{10}=\ds\frac{c_2}{(c_1+c_2)^2},
x_{20}=\ds\frac{c_1}{(c_1+c_2)^2};
\end{equation*}

(ii) The elements of the matrix $A$, which characterize linearized
equation (\ref{31}) in $(x_{10}, x_{20})$ are:
\begin{equation*}\label{35}
\begin{split}
a_{11}=-2k_1c_1(c_1+c_2), a_{12}=-k_1(c_1^2-c_2^2)\\
a_{21}=k_2(c_1^2-c_2^2), a_{22}=-2k_2c_2(c_1+c_2);
\end{split}
\end{equation*}

(iii) The roots of the characteristic equation:
\begin{equation}\label{36}
\mu^2-(a_{11}+a_{22})\mu+a_{11}a_{22}-a_{12}a_{21}=0
\end{equation} have the real part:
\begin{equation*}\label{37}
Re(\mu_{1,2})=-(k_1c_1+k_2c_2)(c_1+c_2);
\end{equation*}

(iv) If $b_{11}=\alpha$, $b_{12}=-\beta$, $b_{21}=\beta$,
$b_{22}=\alpha$, $\beta\neq 0$, then the Lyapunov coefficient of
(SDE) (\ref{3}) is:
\begin{equation}\label{38}
\begin{split}
\lambda\!
&=\!-(k_1c_1+k_2c_2)(c_1+c_2)+\ds\frac{1}{2}(\beta^2-\alpha^2)-
(k_1c_1-k_2c_2)(c_1+c_2)D_2+\\
&+\ds\frac{1}{2}(k_2-k_1)(c_1^2-c_2^2)E_2
\end{split}
\end{equation} where
\begin{equation*}\label{39}
D_2=\int_0^{2\pi}cos(2\theta)p(\theta)d\theta,
E_2=\int_0^{2\pi}sin(2\theta)p(\theta)d\theta
\end{equation*} and
\begin{equation*}\label{40}
\begin{split}
& p(\theta)=kg(\theta),
k=\ds\frac{1}{\int_0^{2\pi}g(\theta)d\theta},\\
&
g(\theta)=\ds\frac{1}{\beta^2}\exp\{\ds\frac{1}{\beta^2}((k_1+k_2)(c_1^2-c_2^2)+\alpha\beta)\theta-
(k_1c_1-k_2c_2)(c_1+c_2)\cos(2\theta)+\\
&+\ds\frac{1}{2}(k_1+k_2)(c_1^2-c_2^2)\sin(2\theta)\}.
\end{split}
\end{equation*}

\end{proposition}

\section{Numerical Simulations.}

\qquad We have done the numerical simulations using a program in
Maple 12. For $c_1=0.2$, $c_2=2$, $k_1=0.2$, $k_2=0.4$, $\beta=2$,
in figure 1 is displayed $(\alpha, \lambda (\alpha))$, where
$\lambda(\alpha)$ is given by (\ref{38}). For $\alpha\in (-\infty,
-1.2) \cup (1.1, \infty)$, the Lyapunov exponent is negative, then
(SDE) has an asymptotically stable stationary state. For $\alpha\in
(-1.2, 1.1)$, the Lyapunov exponent is positive and (SDE) has an
asymptotically unstable stationary state.

\begin{center}\begin{tabular}{ccc}
\\ Fig 1. $(\alpha, \lambda(\alpha))$\\
\includegraphics[width=5cm]{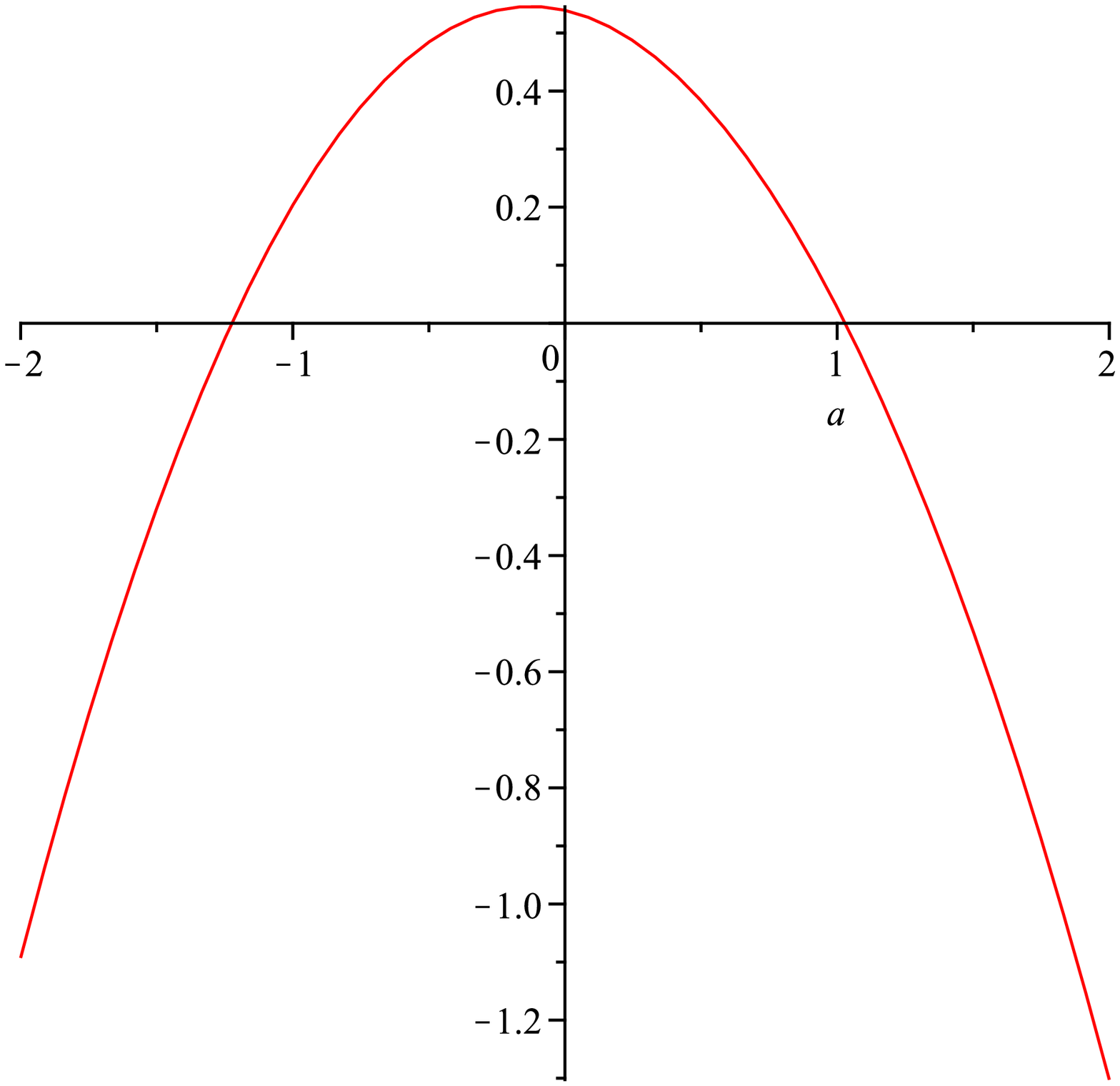}

\end{tabular}
\end{center}

If $\beta$ is a real parameter and $\alpha=2$, in figure 2 we have:
$(\beta, \lambda(\beta))$.

\begin{center}\begin{tabular}{ccc}
\\ Fig 2. $(\beta, \lambda(\beta))$\\
\includegraphics[width=5cm]{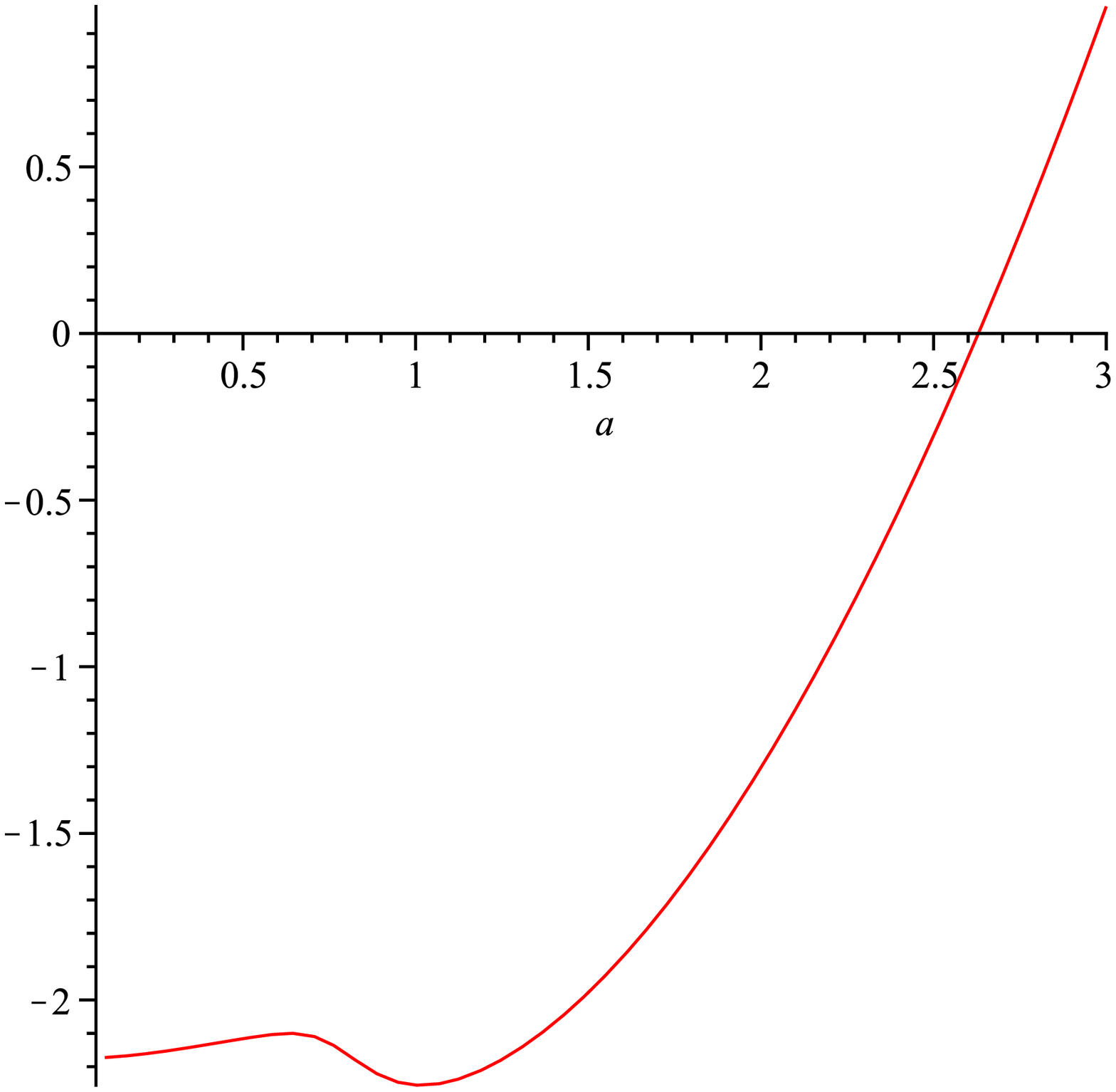}

\end{tabular}
\end{center}

For $\beta \in (-\infty, -2.6)\cup (2.6,\infty)$ the Lyapunov
exponent is positive and (SDE) has an asymptotically unstable
stationary state. For $\beta (-2.6, 2.6)$ the Lyapunov exponent is
negative and (SDE)equation has an asymptotically stable stationary
state.

The Euler second order scheme for (SDE) (\ref{31}) is given by:
\begin{equation*}
\begin{split}
&x_1(n+1)=x_1(n)+h\left (\ds\frac{x_2(n)}{(x_1(n)+x_2(n))^2}-c_1\right )+(b_{11}x_1(n)+b_{12}x_2(n)+\gamma_1)\\
&\cdot
G(n)+b_{11}(b_{11}x_1(n)+b_{12}x_2(n)+\gamma_1)\ds\frac{G(n)^2-h}{2}+
(-\ds\frac{2x_1(n)x_2(n)}{(x_1(n)+x_2(n))^3}\\
&\cdot\left (\ds\frac{x_2(n)}{(x_1(n)+x_2(n))^2}-c_1\right
)+(b_{11}x_1(n)+b_{12}x_2(n)+\gamma_1)\ds\frac{x_1(n)x_2(n)}{(x_1(n)+x_2(n))^3}
)\ds\frac{h^2}{2}+\\
&(b_{11}-\ds\frac{2x_2(n)}{(x_1(n)+x_2(n))^3})(b_{11}x_1(n)+b_{12}x_2(n)+\gamma_1)\ds\frac{hG(n)}{2},
\end{split}
\end{equation*}

\begin{equation*}
\begin{split}
&x_2(n+1)=x_2(n)+h\left (\ds\frac{x_1(n)}{(x_1(n)+x_2(n))^2}-c_2\right )+(b_{21}x_1(n)+b_{22}x_2(n)+\gamma_2)\\
&\cdot
G(n)+b_{22}(b_{21}x_1(n)+b_{22}x_2(n)+\gamma_2)\ds\frac{G(n)^2-h}{2}+
(-\ds\frac{2x_1(n)x_2(n)}{(x_1(n)+x_2(n))^3}\\
&\cdot\left (\ds\frac{x_1(n)}{(x_1(n)+x_2(n))^2}-c_2\right
)+(b_{21}x_1(n)+b_{22}x_2(n)+\gamma_2)\ds\frac{x_1(n)x_2(n)}{(x_1(n)+x_2(n))^3}
)\ds\frac{h^2}{2}+\\
&(b_{21}-\ds\frac{2x_1(n)}{(x_1(n)+x_2(n))^3})(b_{21}x_1(n)+b_{22}x_2(n)+\gamma_2)\ds\frac{hG(n)}{2},
\end{split}
\end{equation*} where $G(n)=w((n+1)h)-w(nh)$, $n=1,2,...$, and $x_i(n)=x_i(nh,\omega)$,
 $i=1,2$.

In figures 3 and 4 the orbits: $(n, x_1(n, \omega))$ for (SDE) and
$(n, x_1(n))$ for (ODE) are displayed:

\begin{center}\begin{tabular}{ccc}
\\ Fig 3. $(n, x_1(n, \omega))$ & Fig 4. $(n, x_1(n))$\\
\includegraphics[width=5cm]{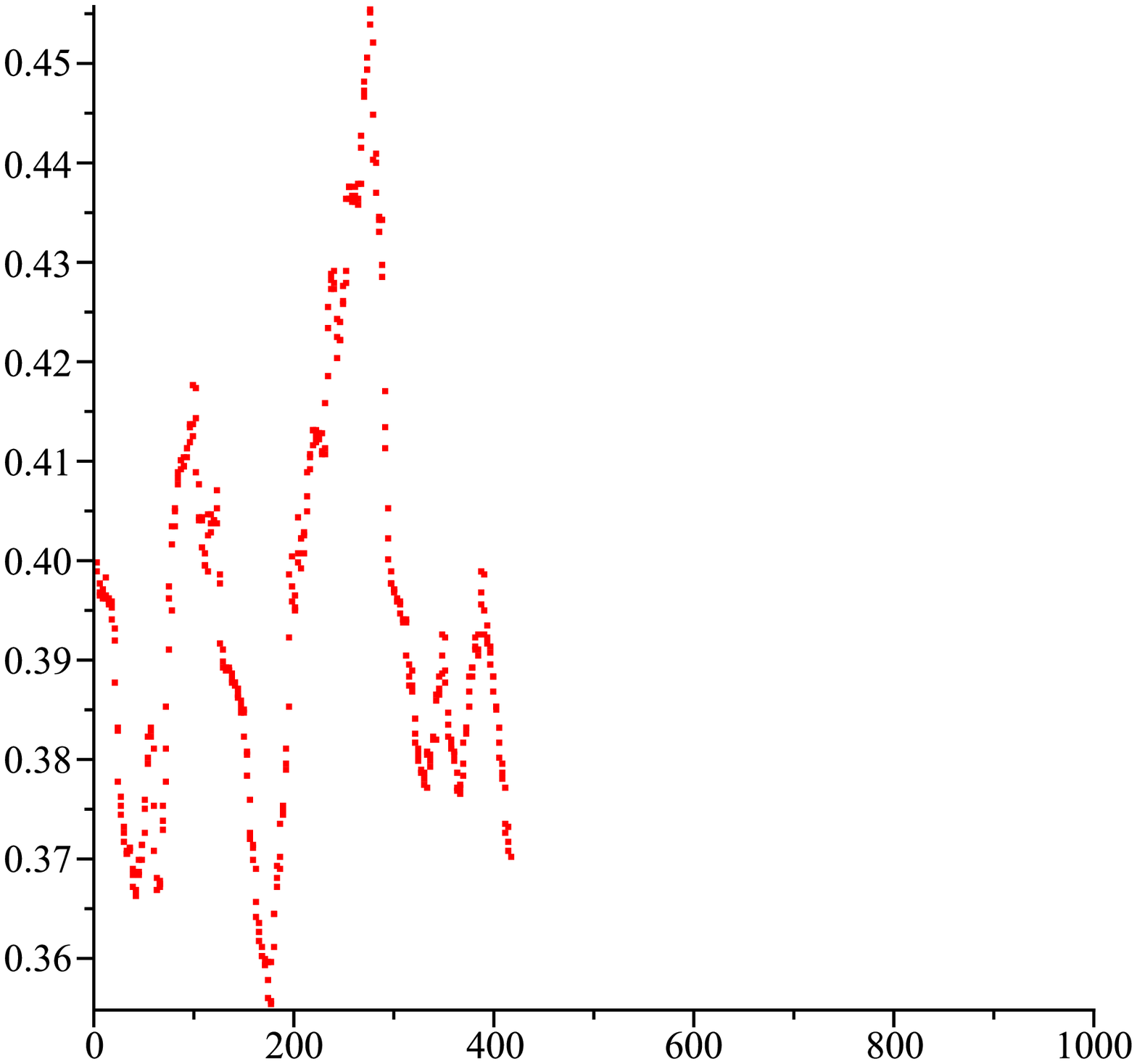} & \includegraphics[width=5cm]{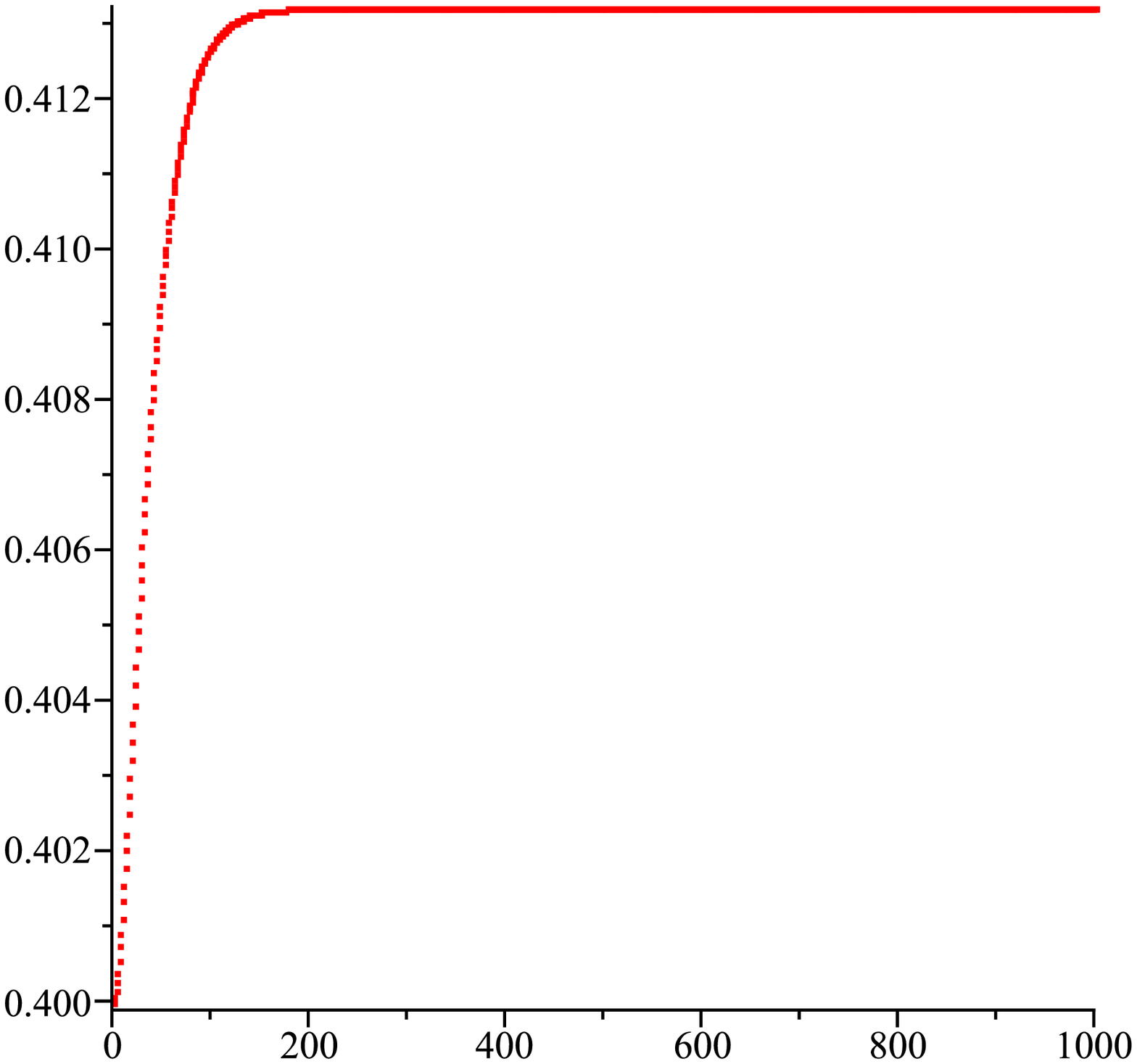}

\end{tabular}
\end{center}

In figures 5 and 6 the orbits: $(n, x_2(n, \omega))$ for (SDE) and
$(n, x_2(n))$ for (ODE) are displayed:

\begin{center}\begin{tabular}{ccc}
\\ Fig 5. $(n, x_2(n, \omega))$ & Fig 6. $(n, x_2(n))$\\
\includegraphics[width=5cm]{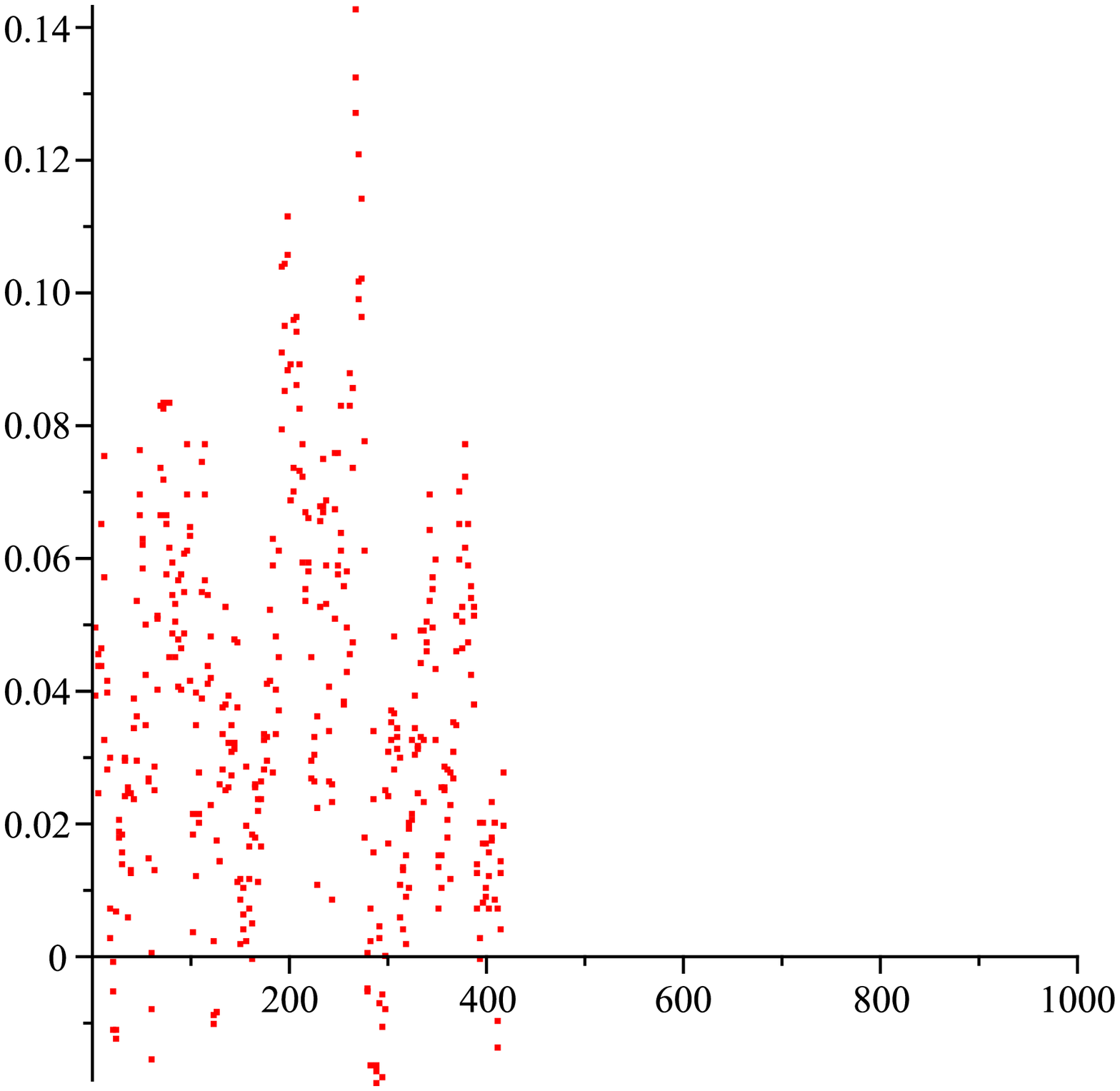} & \includegraphics[width=5cm]{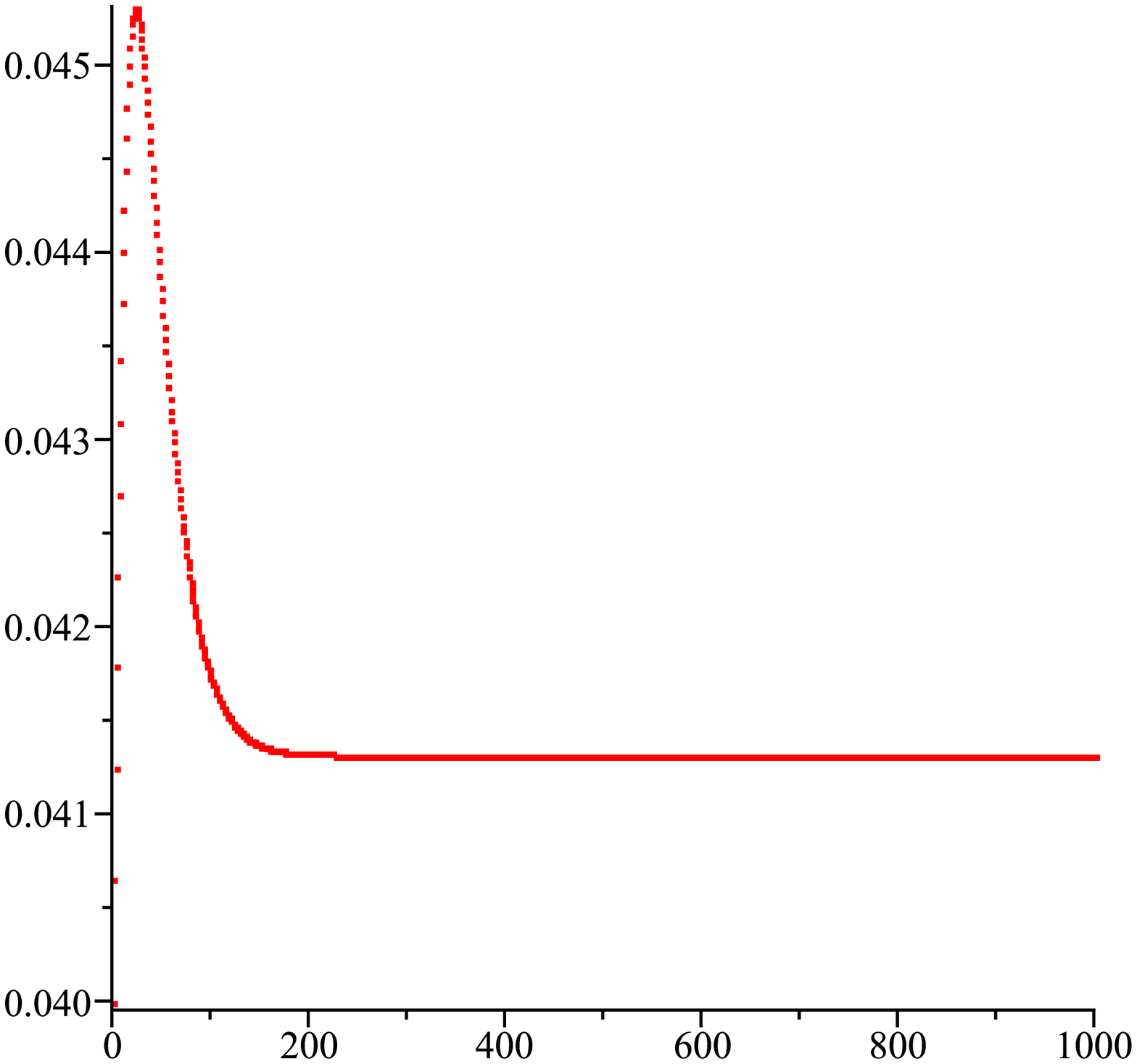}

\end{tabular}
\end{center}

In figures 7 and 8 the orbits: $(x_1(n,\omega), x_2(n, \omega))$ for
(SDE) and  $(x_1(n), x_2(n))$ for (ODE) are displayed:

\begin{center}\begin{tabular}{ccc}
\\ Fig 7. $(x_1(n,\omega), x_2(n,
\omega))$ & Fig 8. $(x_1(n), x_2(n))$\\
\includegraphics[width=5cm]{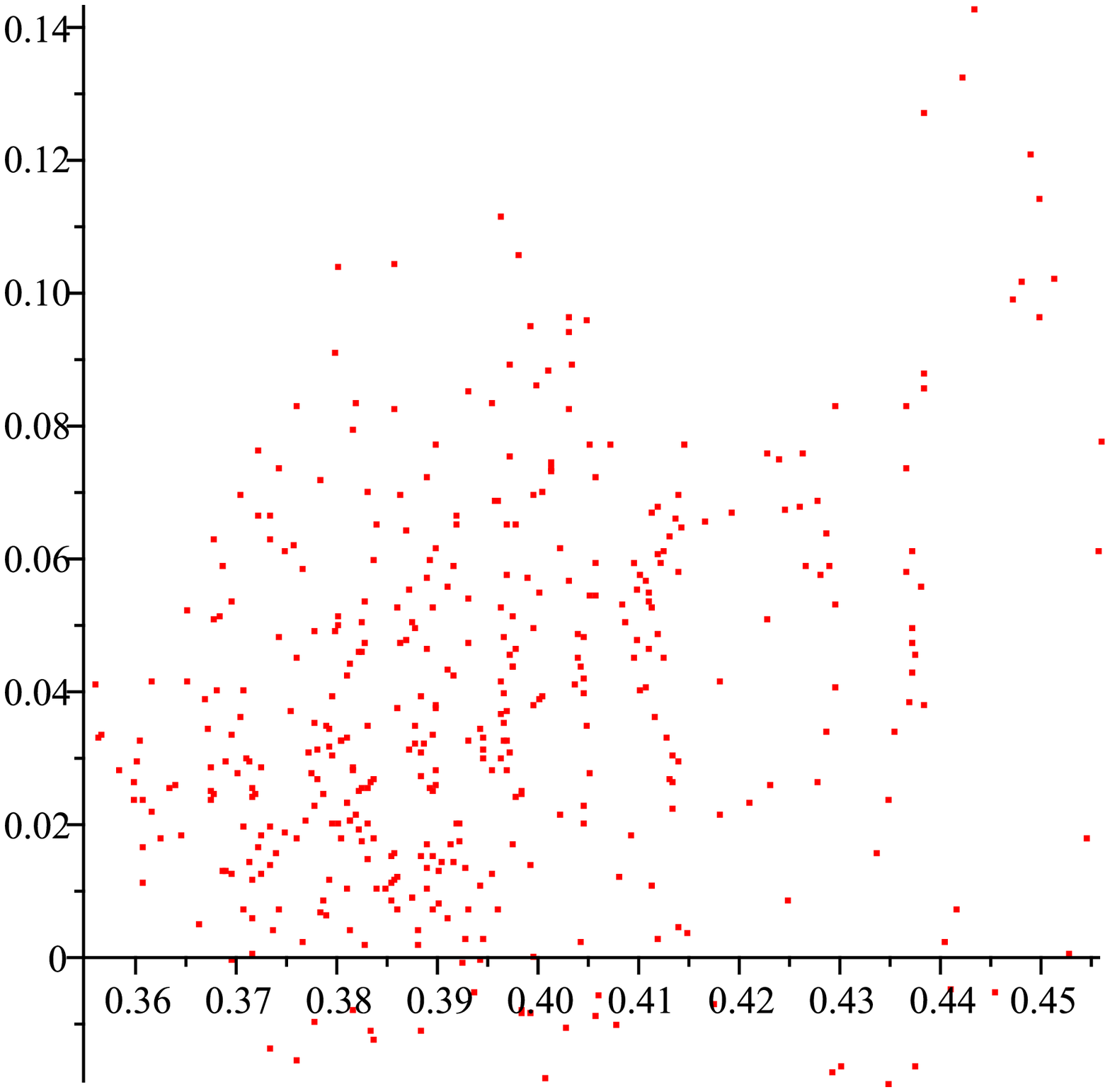} & \includegraphics[width=5cm]{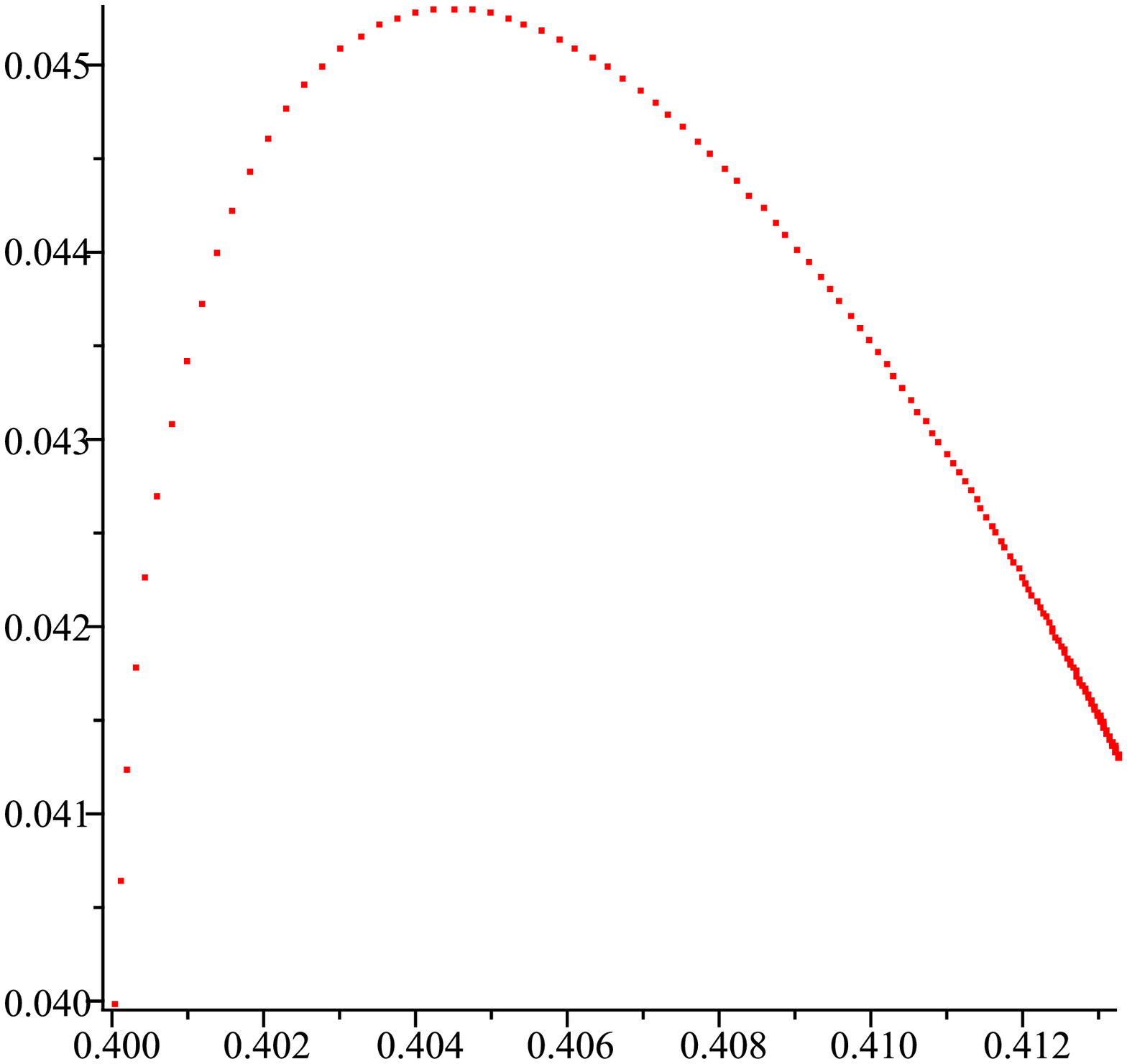}

\end{tabular}
\end{center}

\section{The Lyapunov function method for the stochastic economic game.}
\qquad Theorem \ref{th3} is used for the analysis of the stability
with the help of the Lyapunov function.

Let $V:D=\{[0,\infty)\times \rr ^2\}\to \rr$ be the function given
by:
\begin{equation*}
V(u)=\ds\frac{1}{2}(\omega_1u_1^2+\omega_2u_2^2),
\end{equation*} where $\omega _i>0, i=1,2$.

Using formula (\ref{36f}) for the linearized system of (\ref{31}):
\begin{equation}\label{10f}
\begin{split}
&
du_1(t)=(a_{11}u_1(t)+a_{12}u_2(t))dt+(b_{11}u_1(t)+b_{12}u_2(t))dB(t)\\
&
du_2(t)=(a_{21}u_1(t)+a_{22}u_2(t))dt+(b_{21}u_1(t)+b_{22}u_2(t))dB(t),
\end{split}
\end{equation} we obtain:
\begin{equation}\label{11f}
\begin{split}
LV(u(t))=
& (a_{11}u_1+a_{12}u_2)\omega_1u_1+(a_{21}u_1+a_{22}u_2)\omega_2u_2+\\
&
\ds\frac{1}{2}[(b_{11}u_1+b_{12}u_2)^2\omega_1+(b_{21}u_1+b_{22}u_2)^2\omega_2]=\\
&
=(a_{11}\omega_1+\ds\frac{1}{2}b_{11}^2\omega_1+\ds\frac{1}{2}b_{21}^2\omega_2)u_1^2+(a_{22}\omega_2+\ds\frac{1}{2}b_{12}^2\omega_1+\ds\frac{1}{2}b_{22}^2\omega_2)u_2^2+\\
&
+(a_{12}\omega_1+a_{21}\omega_2+b_{11}b_{12}\omega_1+b_{21}b_{22}\omega_2)u_1u_2.
\end{split}
\end{equation}

If
\begin{equation}\label{12f}
\begin{split}
& A_1=-\ds\frac{a_{21}+b_{21}b_{22}}{a_{12}+b_{11}b_{12}}, \quad
a_{12}+b_{11}b_{12}\neq 0, \omega_1=-A_1\omega_2,\\
& q_1=(a_{11}+ds\frac{1}{2}b_{11})A_1-\ds\frac{1}{2}b_{21}^2,
q_2=-a_{22}-\ds\frac{1}{2}b_{22}^2+\ds\frac{1}{2}b_{12}^2A_1,\\
\end{split}
\end{equation} then from (\ref{11f}) and (\ref{12f}) we get:
$$LV(u)=-q_1\omega_2u_1^2-q_2\omega_2u_2^2.$$

Form the above relations and Theorem \ref{th3} we obtain:

\begin{proposition}
If $b_{ij}, i,j=1,2$ satisfy the relations:
$$a_{12}+b_{11}b_{12}\neq 0, A_1<0, q_1>0, q_2>0,$$ then the trivial
solution of (\ref{10f}) is mean square stable.
\end{proposition}

\section{Conclusions.}

\qquad In the present paper we investigate an economic game with
stochastic dynamics. We focus on a particular game and determine the
Lyapunov exponent for the stochastic system of equations that
describes the mathematical model and the Lyapunov function for the
analysis of the mean square stability. The calculation of the top
Lyapunov exponent allows us to decide whether a stochastic system is
stable or not. Using a program in Maple 12, we display the Lyapunov
exponent and the system orbits. Conditions for the solution of the
stochastic game to be asymptotically mean square stable are
established.

\end{document}